# Logical contradictions in the One-way ANOVA and Tukey-Kramer multiple comparisons tests with more than two groups of observations

Vladimir Gurvich

*National Research University Higher School of Economics, Moscow, Russia*

Mariya Naumova

*Rutgers Business School, Rutgers University, Piscataway, New Jersey, United States*

## Abstract

We show that the One-way ANOVA and Tukey-Kramer ($TK$) tests agree on any sample with two groups. This result is based on a simple identity connecting the Fisher-Snedecor and studentized probabilistic distributions and is proven without any additional assumptions; in particular, the standard ANOVA assumptions (independence, normality, and homoscedasticity ($INAH$)) are not needed. In contrast, it is known that for a sample with $k > 2$ groups of observations, even under the $INAH$ assumptions, with the same significance level $\alpha$, the above two tests may give opposite results:

(i) ANOVA rejects its null hypothesis $H_0^A : \mu_1 = \ldots = \mu_k$, while the $TK$ one, $H_0^{TK}(i,j) : \mu_i = \mu_j$, is not rejected for any pair $i, j \in \{1, \ldots, k\}$;

(ii) the $TK$ test rejects $H_0^{TK}(i,j)$ for a pair $(i,j)$ (with $i \neq j$) while ANOVA does not reject $H_0^A$.

We construct two large infinite pseudo-random families of samples of both types satisfying $INAH$: in case (i) for any $k \geq 3$ and in case (ii) for some larger $k$. Furthermore, in case (ii) ANOVA, being restricted to the pair of groups $(i,j)$, may reject equality $\mu_i = \mu_j$ with the same $\alpha$. This is an obvious contradiction, since $\mu_1 = \ldots = \mu_k$ implies $\mu_i = \mu_j$ for all $i, j \in \{1, \ldots, k\}$. Similar contradictory examples are constructed for the Multivariable Linear

---

*Email addresses:* vgurvich@hse.ru and vladimir.gurvich@gmail.com (Vladimir Gurvich), mnaumova@business.rutgers.edu (Mariya Naumova)

*May 22, 2021*

Regression (MLR). However, for these constructions it seems difficult to verify the Gauss-Markov assumptions, which are standardly required for MLR.

Mathematics Subject Classification: 62 Statistics

---

**1. One-way ANOVA and Tukey-Kramer multiple comparisons tests**

We use standard statistical definitions and notation; the reader can find more details in Miller (1999) or Montgomery (2017).

*1.1. One-way ANOVA*

Consider an arbitrary sample that consists of $k$ groups of randomly chosen real values. A group $j \in \{1, \ldots, k\}$ contains $n_j$ values $x_{\ell j}$ with $\ell = 1, \ldots, n_j$. Then $n = n_1 + \ldots + n_k$ is the total number of values in the sample.

Standardly, $\bar{x}_j$ and $\mu_j$ denote the *sample* and *population means* for $j = 1, \ldots, k$.

We test
$$H_0^A : \mu_1 = \ldots = \mu_k$$
$$H_1^A : \text{not all } \mu_i \text{ are the same}, i = 1, \ldots, k.$$

The *One-way ANOVA* test rejects the null hypothesis $H_0^A$ with significance $\alpha$, that is, with confidence $100(1-\alpha)\%$, if and only if
$$F_{stat} > F_{crit}(\alpha, k-1, n-k),$$
or equivalently, if the $p$-value corresponding to $F_{stat}$ is less than $\alpha$.

Here $F_{crit}(\alpha, k-1, n-k)$ is the critical value of the Fisher-Snedecor distribution corresponding to the significance level $\alpha$, with degrees of freedom of the numerator $df_1 = k - 1$ and of the denumerator $df_2 = n - k$.

The value $F_{stat}$ is given by the ratio
$$F_{stat} = \frac{MS(Tr)}{MSE},$$
where
$$MS(Tr) = \frac{SS(Tr)}{k-1}, \quad SS(Tr) = \sum_{j=1}^{k} n_j (\bar{x}_j - \bar{\bar{x}})^2 = \frac{1}{n} \sum_{j=1}^{k} \sum_{i=j+1}^{k} n_j n_i (\bar{x}_j - \bar{x}_i)^2,$$

$$\bar{\bar{x}} = \frac{1}{n}\sum_{j=1}^{k}\sum_{\ell=1}^{n_j} x_{\ell j} = \frac{1}{n}\sum_{j=1}^{k} n_j \bar{x}_j. \tag{1}$$

$$MSE = \frac{SSE}{n-k}, \quad SSE = \sum_{j=1}^{k}\sum_{i=1}^{n_j}(x_{ij} - \bar{x}_j)^2,$$

Thus, ANOVA rejects $H_0^A$ if and only if

$$MSE < (n(k-1)F_{crit}(\alpha, k-1, n-k))^{-1} \sum_{j=1}^{k}\sum_{i=j+1}^{k} n_j n_i (\bar{x}_j - \bar{x}_i)^2. \tag{2}$$

*1.2. Tukey-Kramer's Test*

For each pair $i, j \in \{1, \ldots k\}$ we test the null hypothesis:

$$H_0^{TK}(i,j) : \mu_i = \mu_j$$
$$H_1^{TK}(i,j) : \mu_i \neq \mu_j$$

for all $i \neq j$.

Tukey (1949) proposed a procedure for testing these hypotheses in case of equal group sizes $n_1 = \ldots = n_k$. Then it was extended in Kramer (1956), Kramer (1957) to arbitrary group sizes. This test, called the Tukey-Kramer ($TK$) test, uses the studentized range statistic

$$Q = \frac{\bar{y}_{max} - \bar{y}_{min}}{\sqrt{\frac{MSE}{n}}},$$

where $\bar{y}_{max}$ and $\bar{y}_{min}$ are the largest and the smallest sample means, out of a collection of $k$ sample means.

$TK$ rejects $H_0^{TK}$ if and only if

$$|\bar{x}_i - \bar{x}_j| > CR(\alpha, k, n, i, j), \tag{3}$$

where the *critical range* ($CR$) is defined by formula

$$CR(\alpha, k, n; i, j) = Q(\alpha, k, n-k)\sqrt{\frac{MSE}{2}\left(\frac{1}{n_i} + \frac{1}{n_j}\right)}, \tag{4}$$

and $Q(\alpha; k, n-k)$ is the critical value of the studentized range $Q$ corresponding to the significance level $\alpha$, with degrees of freedom of numerator $df_1 = k$ and of denominator $df_2 = n - k$.

Equivalently, (3) can be stated as

$$MSE < 2Q^{-2}(\alpha, k, n-k)(\bar{x}_i - \bar{x}_j)^2 \frac{n_i n_j}{n_i + n_j}. \tag{5}$$

*1.3. Comparing ANOVA and Tukey-Kramer tests*

Both ANOVA and *TK* tests are based on the following standard assumptions: independence, normality, and homoscedasticity (*INAH*). The last one means that all groups have equal standard deviations, $\sigma_1 = \ldots = \sigma_k$. Both rejection criteria (2) and (5) are based on these assumptions; see Montgomery (2017), Miller (1999), and Hayter (1984) for more details.

In this note we concentrate on the agreement between the above two tests rather than on their validity. Both inequalities (2) and (5) have the same left-hand side, $MSE$, which can be any number and it is irrelevant for the sake of comparison of two tests.

By definition, $H_0^A$ holds if and only if $H_0^{TK}$ holds for all pairs $(i,j)$ with $i \neq j$. When $k = 2$ there is only one such pair and, hence, the ANOVA and *TK* tests should agree, and indeed they are. In Section 2 we prove it for an arbitrary sample. In particular, even if the *INAH* assumptions are not met, still both tests either reject their null hypothesis or both do not, for any fixed significance level $\alpha$.

However, when $k > 2$, even under *INAH* assumptions, the ANOVA and TK tests may disagree and both cases (i) and (ii) defined in Abstract may take place. Case (i) is not a paradox. Indeed, if $H_0^{TK}(i,j) : \mu_i = \mu_j$ holds with significance slightly larger than $\alpha$ then it is not rejected by the *TK* test. This may hold for all pairs $i,j \in \{1,\ldots,k\}$ with $i \neq j$. Yet, the number of these pairs $\frac{k(k-1)}{2}$ is more than 1 when $k > 2$. So, ANOVA may reject $H_0^A : \mu_1 = \ldots = \mu_k$ with significance $\alpha$.

Somewhat surprisingly, the inverse happens too: $H_0^A$ may hold with a fixed $\alpha$ while $H_0^{TK}$ may be rejected for some pair $(i,j)$ with the same $\alpha$. Such examples are known. Hsu (1996) on page 177 remarks: "An unfortunate common practice is to pursue multiple comparisons only when the null hypothesis of homogeneity is rejected".

We construct two large families of samples of both types considered above. In Section 3.1 we provide two randomly generated samples with three groups in each, $k = 3$, and in Section 3.2 two infinite families of pseudo-random samples with $k \geq 3$ for type (i) and with some larger $k$ for type (ii). It is important to note that these constructions are realized under *INAH* assumptions.

When $k > 2$, formula (4) looks somewhat strange: the critical range is defined for a given pair $(i,j)$ via the value of $MSE$ that depends on all observations, in all $k$ groups. These observations are independent random variables; hence, their values in a group $\ell$ cannot affect the equality $H_0^{TK}(i,j) : \mu_i = \mu_j$ whenever $i,j,\ell$ are pairwise distinct. Moreover, group $\ell$ may be not related to groups $i$ and $j$ at all. Equal or not, mean sales of two New Jersey supermarkets should not depend on the mean accumulation of snow in Siberia. And yet, according to (4), it does.

*1.4. Modified Tukey-Kramer test*

Given significance level $\alpha$, a sample with $k > 2$ groups, and a pair $(i,j)$ with $i \neq j$, let us modify the $TK$ test for $H_0^{TK}(i,j) : \mu_i = \mu_j$ by eliminating all groups but $i$ and $j$ from the sample. Thus, we obtain a new sample with $k' = 2$, $n' = n_i + n_j$, and

$$MSE' = \frac{SSE'}{n-2}, \quad SSE' = \sum_{i'=1}^{n_i}(x_{i'i} - \bar{x}_i)^2 + \sum_{j'=1}^{n_j}(x_{j'j} - \bar{x}_j)^2. \tag{6}$$

Then we define

$$CR'(\alpha, k, n; i, j) = Q(\alpha, 2, n' - 2)\sqrt{\frac{MSE'}{2}\left(\frac{1}{n_i} + \frac{1}{n_j}\right)}. \tag{7}$$

In Section ??, assuming homoscedasticity ($\sigma_1 = \ldots = \sigma_k$) and also that $n_1 = \ldots = n_k$ and $\frac{n}{k}$ is large enough, we will show that $CR' \leq CR$. Hence, the modified $TK$ test rejects $H_0^{TK}(i,j) : \mu_i = \mu_j$ whenever the standard $TK$ test does.

**Remark 1.** *Note that in general, inequality $CR' \leq CR$ may fail since $MSE$ may be much smaller than $MSE'$. Indeed, if $s_\ell = 0$ (resp., small) for all $\ell \in \{1,\ldots,k\}\setminus\{i,j\}$, while $s_i > 0$ and $s_j > 0$ (resp., large) then $SSE' = SSE > 0$ (resp., $SSE - SSE'$ may be an arbitrarily small non-negative number). Notice however that the homoscedasticity assumption might not hold when $s_i$ and $s_j$ differ significantly. Furthermore, $n - k$ may be much larger than $n' - 2$ resulting in $Q(k, n-k)\sqrt{MSE} < Q(2, n'-2)\sqrt{MSE'}$. We leave constructing such examples to the careful reader.*

*1.5. Counter-intuitive examples with symmetric samples of 2 and 3 groups*

Another infinite set of pseudo-random examples will be constructed in Section 4. Given two groups of observations with $n_1 = n_2 = d$, $\bar{x}_1 = -1, \bar{x}_2 = 1$, $\sigma_1 = \sigma_2 = \sigma$, and a significance level $\alpha$, we find $d$, $\sigma$, and $\alpha$ such that ANOVA rejects $H_0 : \mu_1 = \mu_2$ with confidence $1 - \alpha$. Then, we add a third group of observations with $n_3 = c$, $\bar{x}_3 = 0$, $\sigma_3 = \sigma$ and show that $H_0^A : \mu_1 = \mu_2 = \mu_3$ is not rejected by ANOVA with the same confidence when $0 < c < d$.

As we already mentioned, this is a logical contradiction. Let us add that condition $c < d$ looks counter-intuitive. Indeed, $\bar{x}_3 = 0$ and hence group 3 contains values that are typically between $\bar{x}_1 = -1$ and $\bar{x}_2 = 1$, which could be viewed as "an argument" in support of $\mu_1 = \mu_2 = \mu_3$. Furthermore, the larger $c$ is the stronger is this argument.

1.6. *ANOVA is not inclusion monotone on the subsets of its k groups of observations* $\{1, \ldots, k\}$

Given a significance level $\alpha$, a sample with $k > 2$ groups, and a pair $i, j \in \{1, \ldots, k\}$ with $i \neq j$, recall case (ii) (the *TK* test rejects $H_0^{TK}(i,j) : \mu_i = \mu_j$, while ANOVA does not reject $H_0^A : \mu_1 = \ldots = \mu_k$).

Reduce the sample to only two groups $i$ and $j$ eliminating $k-2$ remaining groups and apply to the obtained sample the ANOVA and modified *TK* tests. According to the previous subsection, the latter still rejects the equality $\mu_i = \mu_j$ and, by Theorem 1, ANOVA also rejects it, while $\mu_1 = \ldots = \mu_k$ was not rejected. This is a contradiction.

1.7. *Logical contradictions in F- and t-tests of Multivariable Linear Regression (MLR)*

The general multivariable linear regression model with $k$ predictors $X^1, \ldots, X^k$ and response $Y$ can be written as

$$Y = \beta_0 + \beta_1 X^1 + \ldots + \beta_k X^k + \epsilon.$$

The properties of the estimators of the coefficients $\beta_i$ are derived under the Gauss-Markov Assumptions (GMA), see for example, Wooldridge (2012).

Commonly used tests in regressions analysis are the *F*-test:

$$H_0^F : \beta_1 = \ldots = \beta_k = 0$$
$$H_1^F : \text{at least one } \beta_i \text{ is not } 0, \text{ for } i = 1, \ldots, k,$$

and the *t*-test for individual coefficients $\beta_i$ for $i \in \{1, \ldots, k\}$:

$$H_0^{t_i} : \beta_i = 0$$
$$H_1^{t_i} : \beta_i \neq 0.$$

It is well-known that the *F*- and *t*-tests are equivalent in case of Simple Linear Regression (SLR), that is, when $k = 1$. In this case, the *p*-values of the tests are equal due to identity $F(1, \nu) = t^2(\nu)$ for all natural $\nu$, where $F(1, \nu)$ is a Fisher-Snedecor random variable with $df_1 = 1$ and $df_2 = \nu$, and $t(\nu)$ is a random variable having Student distribution with the degrees of freedom $\nu$.

Yet, for MLR, $k > 1$, logical contradictions similar to ones outlined in Subsections 1.3–1.5 appear. With the same significance level $\alpha$, the *F*- and *t*-tests for MLR may give opposite results:

(j) *F*-test rejects $H_0^F : \beta_1 = \ldots = \beta_k = 0$, while $H_0^{t_i} : \beta_i = 0$ is rejected for no $i \in \{1, \ldots, k\}$;

(jj) $F$-test does not reject $H_0^F$, while $t$-test rejects $H_0^{t_i}$ for some (or even for all) $i \in \{1,\ldots,k\}$.

Similarly to case $(i)$ of ANOVA, case $(j)$ is not a paradox: $H_0^{t_i}$ cannot be rejected with significance $\alpha$ for each particular $i$, but it can be rejected with this significance for at least one $i$. In contrast, case $(jj)$ is an obvious contradiction, since $H_0^F : \beta_1 = \ldots = \beta_k = 0$ implies $H_0^{t_i} : \beta_i = 0$ for every $i \in \{1,\ldots,k\}$.

The corresponding examples are shown in section 5 for $k = 2$ with the following inequalities for the $p$-values:

- $p_{12} > p_1$ and $p_{12} > p_2$ in case $(j)$; see Tables 5A-B in section 5,

- $p_{12} < p_1$ and $p_{12} < p_2$ in case $(jj)$; see Tables 6A-B, 7A-B in section 5,

where $p_{12}$ is the $p$-value of the $F$-test, while $p_1$ and $p_2$ are the $p$-values of the $t$-tests for $\beta_1$ and $\beta_2$, respectively.

Similarly to Subsection 1.5 for ANOVA, we will show that MLR may be not inclusion monotone on the set $\{1,\ldots,k\}$ of its predictors. More precisely, consider MLR $F$-test with $k$ predictors and eliminate $k-1$ of them, all but $i$, getting $k$ SLR problems, one for each predictor $X^i$ and the same response $Y$, for $i \in \{1,\ldots,k\}$. Denote $p'_i$ the $p$-value of the SLR test $i$. (Recall that the $F$- and $t$-tests for SLR are equivalent and have equal $p$-values.)

The example in Table 7A also has the following property: after removing predictor $X^1$, we obtain $p_{12} > 0.05 > p'_2$; see Table 7C. Hence, for significance level $\alpha = 0.05$ the null hypothesis $H_0^F : \beta_1 = \beta_2 = 0$ is not rejected, while $p_2$ and $p'_2$ are both less than 0.05. Thus, Case $(jj)$ holds and, furthermore, both predictors $X^1$ and $X^2$ are not significant, while $X^2$ alone is significant.

Let us remark finaly that our constructions of Sections 1.3–1.5 satisfy INAH assumptions for ANOVA, however, it seems difficult to verify the GMA, which are standardly required for MLR.

## 2. Two groups, $k = 2$

In this case there is a unique pair $(i,j) = (1,2)$ of means and the multiple comparisons turn into a single one. ANOVA and TK tests' null hypotheses $H_0$ coincide stating that $\mu_1 = \mu_2$.

**Theorem 1.** *In case of two groups, k=2, ANOVA and TK tests are equivalent.*

*Proof.* It is enough to show that inequalities (2) and (5) are equivalent when $k = 2$. In this case formulas (2) and (5) can be rewritten as follows:

$$MSE < \frac{n_1 n_2}{n} F_{crit}^{-1}(\alpha, 1, n-2)(\bar{x}_1 - \bar{x}_2)^2,$$

and

$$MSE < \frac{2n_1 n_2}{n} Q^{-2}(\alpha, 2, n-2)(\bar{x}_1 - \bar{x}_2)^2,$$

where

$$MSE = \frac{1}{n-2} \sum_{j=1}^{k} \sum_{i=1}^{n_j} (x_{ij} - \bar{x}_j)^2.$$

Thus, it is sufficient to prove the identity

$$Q^2(\alpha, 2, n-2) = 2F_{crit}(\alpha, 1, n-2),$$

which is implied by the following lemma.

Let $F(1, \nu)$ be a Fisher-Snedecor random variable with $df_1 = 1$ and $df_2 = \nu$, and $Q(2, \nu)$ be a random variable having studentized range distribution with the number of groups $k = 2$ and the degrees of freedom $\nu$.

**Lemma 1.** *Equation* $2F(1, \nu) = Q^2(2, \nu)$ *holds.*

*Proof.* The probability density function of studentized range $Q$ in case $k = 2$ is given by

$$f_Q(q; 2, \nu) = \frac{4\sqrt{2\pi} \left(\frac{\nu}{2}\right)^{\frac{\nu}{2}}}{\Gamma(\frac{\nu}{2})} \int_0^\infty s^\nu \phi(\sqrt{\nu}s) \left[ \int_\infty^\infty \phi(z+qs)\phi(z) dz \right] ds$$

$$= \frac{4\sqrt{2\pi} \left(\frac{\nu}{2}\right)^{\frac{\nu}{2}}}{\Gamma(\frac{\nu}{2})} \int_0^\infty s^\nu \frac{1}{\sqrt{2\pi}} e^{-\frac{\nu s^2}{2}} \left[ \int_\infty^\infty \frac{1}{2\pi} e^{-\frac{(z+qs)^2+z^2}{2}} dz \right] ds;$$

see Tukey (1949). We transform this formula as follows. Substitute $u = \sqrt{2}z$ to obtain

$$f_Q(q; 2, \nu) = \frac{2\sqrt{2} \left(\frac{\nu}{2}\right)^{\frac{\nu}{2}}}{\Gamma(\frac{\nu}{2})} \int_0^\infty s^\nu \frac{1}{\sqrt{2\pi}} e^{-\frac{\nu s^2}{2}} \left[ \int_\infty^\infty \frac{1}{\sqrt{2\pi}} e^{-\frac{\left(u+\frac{\sqrt{2}qs}{2}\right)^2 + \frac{q^2 s^2}{2}}{2}} du \right] ds$$

$$= \frac{2\sqrt{2} \left(\frac{\nu}{2}\right)^{\frac{\nu}{2}}}{\sqrt{\pi}\Gamma(\frac{\nu}{2})} \int_0^\infty s^\nu e^{-\frac{\left(\nu+\frac{q^2}{2}\right)s^2}{2}} ds.$$

Then by substitution $t = s^2$,

$$f_Q(q; 2, \nu) = \frac{\left(\frac{\nu}{2}\right)^{\frac{\nu}{2}}}{\sqrt{\pi}\Gamma(\frac{\nu}{2})} \int_0^\infty t^{\frac{\nu-1}{2}} e^{-\frac{\left(\nu + \frac{q^2}{2}\right)t}{2}} dt,$$

and by substitution $y = \frac{\left(\nu + \frac{q^2}{2}\right)t}{2}$,

$$\begin{aligned}
f_Q(q; 2, \nu) &= \frac{\sqrt{2}\nu^{\frac{\nu}{2}}}{\sqrt{\pi}\Gamma(\frac{\nu}{2})\left(\nu + \frac{q^2}{2}\right)^{\frac{\nu+1}{2}}} \int_0^\infty y^{\frac{\nu+1}{2}-1} e^{-y} dy \\
&= \frac{\sqrt{2}\nu^{\frac{\nu}{2}}}{\sqrt{\pi}\Gamma(\frac{\nu}{2})\left(\nu + \frac{q^2}{2}\right)^{\frac{\nu+1}{2}}} \Gamma\left(\frac{\nu+1}{2}\right).
\end{aligned} \quad (8)$$

For $X = \frac{Q^2}{2}$, we obtain

$$P[X \leq x] = P\left[\frac{Q^2}{2} \leq x\right] = P\left[Q \leq \sqrt{2x}\right],$$

and then

$$f_X(x) = \frac{d}{dx} P\left[Q \leq \sqrt{2x}\right] = \frac{1}{\sqrt{2x}} f_Q(\sqrt{2x}; 2, \nu),$$

which by (8) implies that

$$\begin{aligned}
f_X(x) &= \frac{1}{\sqrt{2x}} \frac{\sqrt{2}\nu^{\frac{\nu}{2}}}{\sqrt{\pi}\Gamma(\frac{\nu}{2})(\nu + x)^{\frac{\nu+1}{2}}} \Gamma\left(\frac{\nu+1}{2}\right) \\
&= \frac{1}{B(\frac{1}{2}, \frac{\nu}{2})} (\nu x)^{-\frac{1}{2}} \left(1 + \frac{x}{\nu}\right)^{-\frac{\nu+1}{2}},
\end{aligned} \quad (9)$$

where $B(a, b)$ is the beta function.

It is well-known (see for example Fisher (1992)) that (9) defines the probability density function of the Fisher-Snedecor distribution with degrees of freedom of the numerator $df_1 = 1$ and of the denumerator $df_2 = \nu$. □

This proves the Theorem. □

Note that Theorem 1 holds for an arbitrary sample. In particular, the $p$-values for ANOVA and $TK$ tests are equal regardless the validity of assumptions *INAH*.

## 3. Some cases when ANOVA and *TK* tests disagree

In this section we provide several examples where the considered two tests disagree: (i) $H_0^A$ is rejected while $H_0^{TK}$ is not, or (ii) vice versa. In Section 3.1 we provide two randomly generated samples illustrating (i) and (ii) with three groups in each, $k = 3$; and in Sections 3.2 and 3.3 we construct two infinite families of pseudo-random samples with $k \geq 3$ for (i) and with some larger $k$ for (ii).

### 3.1. Two examples with 3 groups

Using R, we generated two random samples with $k = 3$ groups, $n_1 = n_2 = n_3 = 10$, from Normal distributions with parameters $\mu_1 = 10$, $\mu_2 = 25$, $\mu_3 = 40$, and $\sigma_1 = \sigma_2 = \sigma_3 = 25$.

*Case 1: ANOVA rejects $H_0^A$ while TK does not reject $H_0^{TK}$*

Table 1A: Generated random sample for Case 1

| Group 1 | Group 2 | Group 3 |
|---|---|---|
| 33.73617429 | 41.34327861 | 1.949654854 |
| 6.532109599 | -2.29596015 | 64.73534452 |
| -15.87068125 | 37.80911436 | 17.47791461 |
| 24.41853292 | -38.5488504 | 43.91077426 |
| 32.52469512 | 28.81447508 | 15.70006485 |
| -36.67775074 | 91.99464773 | 54.51355702 |
| 4.946144821 | 54.96895462 | 31.54941908 |
| -11.48789077 | 77.16877 | 0.819316511 |
| 32.04750431 | 95.1318948 | 73.84330239 |
| -24.02530782 | 6.722405014 | 45.92357859 |

Table 1B: ANOVA Table for the example in Case 1

|  | Df | Sum Sq | Mean Sq | F value | Pr(>F) |
|---|---|---|---|---|---|
| group | 2 | 7159.763 | 3579.881 | 3.39789 | 0.04828 |
| Residuals | 27 | 28446.164 | 1053.562 |  |  |

Table 1C: The results of $TK$ test for the example in Case 1

| group | diff | lwr | upr | p adj |
|---|---|---|---|---|
| Group2-Group1 | 34.696519918 | -1.294542536 | 70.68758237 | 0.0604457982 |
| Group3-Group1 | 30.427939620 | -5.563122834 | 66.41900207 | 0.1095211230 |
| Group3-Group2 | -4.268580298 | -40.259642752 | 31.72248216 | 0.9535328433 |

Let us fix the significance level $\alpha = 0.05$, then $TK$ does not distinguish any pair $\mu_i$ and $\mu_j$, while Table 1B shows that ANOVA test rejects the null hypothesis $H_0^A$.

*Case 2: ANOVA does not reject $H_0^A$ while $TK$ rejects $H_0^{TK}$*

Table 2A: Generated random sample for Case 2

| Group 1 | Group 2 | Group 3 |
|---|---|---|
| 19.65656273 | 30.47282693 | 97.66594506 |
| 31.63471018 | 2.359493274 | 37.29706457 |
| 5.474716521 | 25.94822801 | 37.28238885 |
| 7.325738946 | -6.706730014 | -9.215515132 |
| 47.16633 | 56.00337827 | 44.75306142 |
| -28.99487682 | 22.37945513 | 72.60365833 |
| 14.99564807 | 73.81358543 | 21.39501942 |
| 48.12035772 | 5.44699726 | 71.5651277 |
| 25.54178184 | -3.745973145 | 63.33149261 |
| -16.61305101 | 48.61987107 | 26.01262136 |

Table 2B: ANOVA Table for the example in Case 2

|  | Df | Sum Sq | Mean Sq | F value | Pr(>F) |
|---|---|---|---|---|---|
| group | 2 | 4948.742 | 2474.371 | 3.20804 | 0.056236 |
| Residuals | 27 | 20825.208 | 771.304 |  |  |

Table 2C: The results of $TK$ test for the example in Case 2

| group | diff | lwr | upr | p adj |
|---|---|---|---|---|
| Group2-Group1 | -10.0283214 | -40.8231285393 | 20.76648573 | 0.7017162217 |
| Group3-Group1 | 30.8382946 | 0.0434874678 | 61.63310174 | 0.0496248072 |
| Group3-Group2 | 20.8099732 | -9.9848339369 | 51.60478033 | 0.2327921189 |

Let us fix again the significance level $\alpha = 0.05$. Then Table 2C shows that $TK$ distinguishes $\mu_1$ and $\mu_3$ at significance level $\alpha = 0.05$. In contrast, Table 2B shows that ANOVA test does not reject $H_0^A$ for the same $\alpha$. In this case we can apply the approach suggested in Sections 1.4–1.5. Let us reduce the sample by eliminating group 2 and apply the ANOVA and (modified) $TK$ test.

Table 2B': ANOVA Table for groups 1 and 3 in Case 2

|  | Df | Sum Sq | Mean Sq | F value | Pr(>F) |
|---|---|---|---|---|---|
| group | 1 | 4755.002 | 4755.002 | 6.044 | 0.024315 |
| Residuals | 18 | 14161.164 | 786.731 |  |  |

Table 2C': The results of modified $TK$ test for groups 1 and 3 in Case 2

| group | diff | lwr | upr | p adj |
|---|---|---|---|---|
| Group3-Group1 | 30.8382946 | 4.484804399 | 57.1917848 | 0.024314834 |

By Theorem 1, these two tests are equivalent: $p$-value is 0.024315 (see Tables 2B', 2C') for the equality $\mu_1 = \mu_2$ in both cases. Yet, for the original sample of 3 groups $p$-value was 0.056236 for the equality $\mu_1 = \mu_2 = \mu_3$. This is an obvious contradiction: ANOVA rejects $\mu_1 = \mu_2$ with confidence 97.5% but cannot reject the stronger statement $\mu_1 = \mu_2 = \mu_3$ (which is easier to do) even with confidence 95%.

Recall that this example was generated by R under *INAH* assumptions. This did not take too many trials: with given parameters $k = 3$, $n_1 = n_2 = n_3 = 10$, $\mu_1 = 10, \mu_2 = 25, \mu_3 = 40$, and $\sigma_1 = \sigma_2 = \sigma_3 = 25$, about each 20 trials provide an example with such properties.

*3.2. Two large families of examples with k groups*

In this subsection we consider samples with $k$ groups such that $n$ is divisible by $k$, and

$$n_1 = \ldots = n_k = \frac{n}{k}, \tag{10}$$

$$s_1 = \ldots = s_k = s, \tag{11}$$

where $s_i$ is the standard deviation of the $i$th group. In this case we have

$$SSE = k\left(\frac{n}{k} - 1\right)s^2, \quad MSE = \frac{SSE}{n - k} = s^2. \tag{12}$$

*Case 1: ANOVA rejects $H_0^A$ while TK does not reject $H_0^{TK}$*

By (2) and (5), this happens if and only if

$$(n(k-1)F_{crit}(\alpha, k-1, n-k))^{-1} \sum_{j=1}^{k} \sum_{i=j+1}^{k} n_j n_i (\bar{x}_j - \bar{x}_i)^2 \qquad (13)$$
$$> MSE > 2Q^{-2}(\alpha, k, n-k)(\bar{x}_i - \bar{x}_j)^2 \frac{n_i n_j}{n_i + n_j},$$

which implies

$$Q^2(\alpha, k, n-k) \sum_{j=1}^{k} \sum_{i=j+1}^{k} n_j n_i (\bar{x}_j - \bar{x}_i)^2 \qquad (14)$$
$$> 2n(k-1)F_{crit}(\alpha, k-1, n-k)(\bar{x}_i - \bar{x}_j)^2 \frac{n_i n_j}{n_i + n_j}.$$

Consider any sample with $k$ groups, $k$ is even, satisfyng (10), (11), and

$$\bar{x}_1 = \ldots = \bar{x}_{\frac{k}{2}} = 1, \bar{x}_{\frac{k}{2}+1} = \ldots = \bar{x}_k = 0. \qquad (15)$$

By (12), $MSE = s^2$, and (14) can be rewritten as

$$\left(\frac{n}{k}\right)^2 \frac{k^2}{4} Q^2(\alpha, k, n-k) > 2 \frac{n^3}{k^2 \left(2\frac{n}{k}\right)} F_{crit}(\alpha, k-1, n-k),$$

which can be simplified to

$$G(\alpha, k, n-k) = \frac{Q^2(\alpha, k, n-k)}{F_{crit}(\alpha, k-1, n-k)} - 4\left(1 - \frac{1}{k}\right) > 0. \qquad (16)$$

Function $G(\alpha, k, n-k)$ has the following properties:

1. $G(\alpha, k, n-k) \equiv 0$ if $k = 2$;

2. $G(\alpha, k, n-k)$ is monotone increasing with respect to $n-k$ and converging as $n \to \infty$ for each $k$;

3. $G(\alpha, k, n-k) > 0$ for $k \geq 3$ and all $n - k \geq 0$ for $\alpha = 0.005, 0.01, 0.025, 0.05, 0.1, 0.25,$ and $0.5$.

It is not our goal to study function $G(\alpha, k, n-k)$ in detail, we are primarily interested only in its positivity, required by condition (16). The required inequality (16) holds for any $k \geq 3$.

Given an even $k$ and $n$ divisible by $k$, we generate a desired pseudo-random sample as follows. It satisfies (10), (11), (15), and in addition, whenever (16) holds, we still can choose $s^2 = MSE$ satisfying (13).

*Case 2: ANOVA does not reject $H_0^A$ while TK rejects $H_0^{TK}$*

By (2) and (5), this happens if and only if

$$(n(k-1)F_{crit}(\alpha, k-1, n-k))^{-1} \sum_{j=1}^{k} \sum_{i=j+1}^{k} n_j n_i (\bar{x}_j - \bar{x}_i)^2 \qquad (17)$$
$$< MSE < 2Q^{-2}(\alpha, k, n-k)(\bar{x}_i - \bar{x}_j)^2 \frac{n_i n_j}{n_i + n_j},$$

which implies

$$Q^2(\alpha, k, n-k) \sum_{j=1}^{k} \sum_{i=j+1}^{k} n_j n_i (\bar{x}_j - \bar{x}_i)^2 \qquad (18)$$
$$< 2n(k-1)F_{crit}(\alpha, k-1, n-k)(\bar{x}_i - \bar{x}_j)^2 \frac{n_i n_j}{n_i + n_j}.$$

Note that if $k = 2$ we obtain (8).

Consider any sample with $k$ groups satisfying (10), (11), and

$$\bar{x}_1 = 1, \bar{x}_2 = \ldots = \bar{x}_k = 0. \qquad (19)$$

Then (18) turns into

$$\left(\frac{n}{k}\right)^2 (k-1) Q^2(\alpha, k, n-k) < \frac{2n(k-1) \left(\frac{n}{k}\right)^2}{2\frac{n}{k}} F_{crit}(\alpha, k-1, n-k),$$

which simplifies to

$$H(\alpha, k, n-k) = \frac{Q^2(\alpha, k, n-k)}{F_{crit}(\alpha, k-1, n-k)} - k < 0. \qquad (20)$$

Since $s$ can be chosen arbitrarily, we can always find $MSE$ satisfying (17) whenever (18) holds.

Function $H(\alpha, k, n-k)$ shares properties (j), (jj) of $G(\alpha, k, n-k)$, and $H(\alpha, k, n-k) > 0$ for sufficiently small $k$, and $H(\alpha, k, n-k) < 0$ for sufficiently large $k$. Again, it is not our goal to study $H(\alpha, k, n-k)$ in detail since we are primarily interested only in its negativity required by condition (20).

The signs of $H(\alpha, k, n-k)$ depending on $k$ are shown in Table 3 for some values of $\alpha$. The second (resp., third) column contains the values of $k$ such that $H(\alpha, k, n-k) > 0$ (resp., $H(\alpha, k, n-k) < 0$) for all $n$. Missing values of $k$ correspond to the cases when the sign of $H(\alpha, k, n-k)$ depends on $n$.

Table 3: The signs of $H(\alpha, k, n-k)$ for selected $\alpha$ depending on $k$

| $\alpha$ | $H(\alpha, k, n-k) > 0$ | $H(\alpha, k, n-k) < 0$ |
|---|---|---|
| 0.005 | $3 \leq k \leq 10$ | $k \geq 14$ |
| 0.01 | $3 \leq k \leq 10$ | $k \geq 14$ |
| 0.025 | $3 \leq k \leq 10$ | $k \geq 13$ |
| 0.05 | $3 \leq k \leq 10$ | $k \geq 12$ |
| 0.1 | $3 \leq k \leq 10$ | $k \geq 12$ |
| 0.25 | $3 \leq k \leq 10$ | $k \geq 11$ |
| 0.5 | $3 \leq k \leq 9$ | $k \geq 10$ |

One can see that the required inequality (20) holds when the number of groups $k$ is large enough.

Given $k$ and $n$ divisible by $k$, we generate a desired preudo-random sample as follows. It satisfies (10), (11), (19), and in addition, whenever (20) holds, we still can choose $s^2 = MSE$ satisfying (17).

**Remark 2.** *We variate the choice of sample means in (15) and (19) to increase the feasible area for (16) and (20), respectively. Obviously, $k \geq 4\left(1 - \frac{1}{k}\right)$ and equality holds if and only if $k = 2$.*

**Remark 3.** *We can extend considerably the family of the constructed examples by relaxing equalities (11), (15), and (19), and replacing them by approximate equalities.*

*3.3. Critical range in Modified TK test*

In Section 1.4 we modified the standard $TK$ multiple comparisons test replacing it by the pairwise comparison version as follows. Given significance level $\alpha$, a sample with $k > 2$ groups, and a pair $(i, j) \in \{1, \ldots, k\}$ with $i \neq j$, consider the null hypothesis for the corresponding two groups, $H_0^{TK}(i, j) : \mu_i = \mu_j$ and eliminate all groups but $i$ and $j$ from the sample, obtaining a new one with $k' = 2$, $n' = n_i + n_j$. For the standard and modified $TK$ tests, the critical ranges $CR = CR(\alpha, k, n-k; i, j)$ and $CR' = CR'(\alpha, k, n-k; i, j)$ and the corresponding values of $MSE$ and $MSE'$ are given by formulas (4), (7), (1), and (6).

We are looking for conditions implying the inequality $CR' \leq CR$, in which case the modified $TK$ test rejects $H_0^{TK}(i, j)$ whenever the standard one does. In general this inequality may fail; see Remark 1.

Let us assume $INAH$, and in addition (10), (11). As we know, in this case $MSE = MSE' = s^2$ and formulas for $CR$ and $CR'$ are simplified as follows:

Table 4: Conditions for monotone increasing of the studentized range $Q(\alpha, 2\ell, \nu\ell)$

| $\alpha$ | 0.005 | 0.01 | 0.025 | 0.05 | 0.1 | 0.25 | 0.5 |
|---|---|---|---|---|---|---|---|
| | $\nu \geq 7$ | $\nu \geq 5$ | $\nu \geq 4$ | $\nu \geq 3$ | $\nu \geq 2$ | $\nu \geq 1$ | $\nu \geq 1$ |

$$CR = Q(\alpha, k, n-k)s\sqrt{\frac{k}{n}},$$

$$CR' = Q(\alpha, 2, n'-k')s\sqrt{\frac{k'}{n'}}$$

$$= Q\left(\alpha, 2, 2\frac{n}{k} - 2\right) s\sqrt{\frac{2}{2n/k}} = Q\left(\alpha, \frac{k}{k/2}, \frac{n-k}{k/2}\right) s\sqrt{\frac{k}{n}}.$$

Thus, in the considered case,

$$\frac{CR'}{CR} = \frac{Q(\alpha, 2, \nu)}{Q(\alpha, 2\ell, \nu\ell)},$$

where $\ell = \frac{k}{2} \geq 1$ and $\nu = n' - k' = 2\left(\frac{n}{k} - 1\right)$.

The critical value of the studentized range $Q(\alpha, 2\ell, \nu\ell)$ monotone increases with $\ell$ when $\nu = 2\left(\frac{n}{k} - 1\right)$ is large enough; see Table 4.

In these cases, $CR' \leq CR$ and, hence, conclusions of Section 1.5 are applicable. Recall the construction of Section 3.2 Case 2, in which ANOVA does not reject $H_0^A : \mu_1 = \ldots = \mu_k$, while the standard $TK$ test rejects $H_0^{TK}(i,j) : \mu_i = \mu_j$. This pseudo-random construction satisfies $INAH$. Let us apply ANOVA and $TK$ tests to the reduced sample that consists of only two groups $i$ and $j$, with the remaining $k - 2$ groups eliminated. By the above arguments, the modified $TK$ test still rejects its hypothesis $H_0^{TK}(i,j) : \mu_i = \mu_j$ and, by Theorem 1, ANOVA rejects it too. However, ANOVA does not reject a stronger hypothesis $H_0^A : \mu_1 = \ldots = \mu_k$, with the same significance level $\alpha$, which is an obvious contradiction.

## 4. Symmetric samples with 2 and 3 groups

### 4.1. Two groups

Consider two groups 1 and 2 with $d$ observations in each, that is, $k = 2$, $n_1 = n_2 = d$, $n = n_1 + n_2 = 2d$, with means $\bar{x}_1 = -1$, $\bar{x}_2 = -1$, and standard deviations $\sigma_1 = \sigma_2 = \sigma$. We can assume that $INAH$ holds.

Obviously, $SS(Tr) = MSR = 2d$; furthermore, by (1),

$$SSE = \sigma^2, \ MSE = SSE/(n-k) = \sigma^2/(2m-2),$$

$$F_{stat} = MSR/MSE = 4m(m-1)\sigma^{-2}.$$

By (2), ANOVA rejects its null hypothesis $H_0^A : \mu_1 = \mu_2$ if and only if

$$4m(m-1)\sigma^{-2} > F_{crit}(\alpha, k-1, n-k) = F_{crit}(\alpha, 1, 2d-2) \qquad (21)$$

As for the $TK$ (in this case just Tukey test), we have $|\bar{x}_1 - \bar{x}_2| = 2$, and by (4) the critical range is given by formula

$$CR(1,2) = Q(\alpha, 2, 2d-2)\sqrt{MSE/d} = Q(\alpha, 2, 2d-2)\sigma/\sqrt{2m(m-1)}.$$

Thus, by (3), the Tukey test rejects $\mu_1 = \mu_2$ if and only if

$$8m(m-1)\sigma^{-2} > Q^2(\alpha, 2, 2m-2). \qquad (22)$$

Criteria (21) and (22) are equivalent, by Lemma 1.

### 4.2. Three groups

Let us add to groups 1 and 2 one more group 3 of $c$ observations getting $k = 3$ and $n = n_1 + n_2 + n_3 = 2d + c$. Furthermore, set $\bar{x}_3 = 0$ and $\sigma_3 = \sigma$, and assume that $INAH$ holds.

Then, by (1) and (2),

$$SS(Tr) = 2d, \ MSR = SSR/(k-1) = d;$$

$$SSE = \sigma^2, \ MSE = SSE/(n-k) = \sigma^2/(2d+c-3),$$

$$F_{stat} = MSR/MSE = d(2d+c-3)\sigma^{-2}.$$

Thus, ANOVA rejects its null hypothesis $H_0^A : \mu_1 = \mu_2 = \mu_3$ if and only if

$$d(2d+c-3)\sigma^{-2} > F_{crit}(\alpha, k-1, n-k) = F_{crit}(\alpha, 2, 2d+c-3). \qquad (23)$$

As for the $TK$ test, we have $|\bar{x}_1 - \bar{x}_2| = 2$, while by (4), the critical range

$$CR(1,2) = Q(\alpha, k, n-k)\sqrt{MSE/d} = \sigma Q(\alpha, 3, 2d+c-3)/\sqrt{d(2d+c-3)}.$$

Thus, by (3), the $TK$ test rejects $\mu_1 = \mu_2$ if and only if

$$4d(2d+c-3)\sigma^{-2} > Q^2(\alpha, 3, 2d+c-3). \qquad (24)$$

### 4.3. ANOVA for $k = 2$ and $k = 3$

ANOVA rejects $\mu_1 = \mu_2$ and does not reject $\mu_1 = \mu_2 = \mu_3$ if and only if (21) holds while (23) fails. It may happen, for some $\sigma$, if and only if

$$\frac{F_{crit}(\alpha, 1, 2a)}{F_{crit}(\alpha, 2, 2a+b)} < 2\frac{2a}{2a+b}, \tag{25}$$

where $a = d - 1$, $b = c - 1$. Obviously, the set of feasible $\sigma$ is an interval.

Consider $\alpha = 0.05$. Inequality (25) holds whenever $b < a$, or equivalently, $c < d$. It seems that (25) can be solved, with respect to $a$ and $b$, explicitly. Consider the following sequence of positive integers $S = (6, A^2, B^4, A, B^5, (A, B^6)^\infty)$, where $A = (8, 9, 8, 8, 9)$, $B = (8, 9, 8, 8, 9, 8, 8, 9)$; a power denotes the number of repetitions. Thus, $S$ is a quasi-periodic sequence with the period $(A, B^6)$ of length $5 + 8 \cdot 6 = 53$. To each $a$ we assign a nonnegative integer $a(s)$ uniquely defined by the inequalities

$$\sum_{i=1}^{a(s)} s_i \leq a < \sum_{i=1}^{a(s)+1} s_i.$$

Then, (25) holds if and only if $b < a + a(s)$. This criterion is confirmed by computations for $a \leq 500$. We conjecture that it holds for all $a$ and that similar criteria hold for all $\alpha$.

### 4.4. TK test for $k = 2$ and $k = 3$

$TK$ rejects equality $\mu_1 = \mu_2$ for $k = 2$ and does not reject it for $k = 3$ if and only if (22) holds while (24) fails. This holds, for some $\sigma$, if and only if

$$\frac{Q^2(\alpha, 3, 2a+b)}{Q^2(\alpha, 2, 2a)} < \frac{2a+b}{2a}, \tag{26}$$

Consider $\alpha = 0.05$. Inequality (26) holds whenever $b \geq a$, or equivalently, $c \geq d$.

Again, it seems that (26) can be solved, with respect to $a$ and $b$, explicitly. Consider sequence of positive integers $S = (3, 7, (8, 7^7, 8, 7^6)^\infty)$. It is quasi-periodic with the period $(8, 7^7, 8, 7^6)$ of length 15.

Then, (26) holds if and only if $b \geq a - a(s)$. This criterion is confirmed by computations for $a \leq 500$. We conjecture that it holds for all $a$ and that similar criteria hold for all $\alpha$.

## 5. Logical contradictions in Multivariable Linear Regression

Here we provide examples announced in Section 1.7.

*Construction for case (j): F-test rejects $H_0^F : \beta_1 = \beta_2 = 0$, while $H_0^{t_1} : \beta_1 = 0$, $H_0^{t_2} : \beta_2 = 0$ are not rejected by t-tests with the same significance $\alpha = 0.05$*

Table 5A: Generated random sample for Case (j)

| X1 | X2 | Y |
|---|---|---|
| 1.713673333 | 0.891652019 | 1.718488057 |
| 0.932830925 | 0.353231823 | 1.311861467 |
| -0.053673724 | 1.132586717 | 1.903344806 |
| 1.055482137 | 0.248411619 | 1.582305067 |
| -0.248355435 | -0.174256727 | 2.607296494 |
| -0.004449867 | 0.115550588 | 2.352411276 |
| 0.086988258 | -0.833496007 | 2.558602277 |
| 0.687284914 | -0.417171685 | 1.721811264 |
| -0.253474712 | 0.045371123 | 1.982673543 |
| 0.135747949 | -0.145817805 | 2.309533234 |

Table 5B: Regression output for the sample in Table 5A

| Coefficients | Estimate | Std. Error | t value | $Pr(>|t|)$ |
|---|---|---|---|---|
| (Intercept) | 2.193398826 | 0.12112005 | 18.10929596 | 3.87188E-07 |
| X1 | -0.397578417 | 0.173782994 | -2.28778667 | 0.055989888 |
| X2 | -0.225853392 | 0.196597153 | -1.148813135 | 0.288369425 |

Residual standard error: 0.321371159 on 7 degrees of freedom
Multiple R-squared: 0.578426902, Adjusted R-squared: 0.457977446
F-statistic: 4.802237548 on 2 and 7 DF, p-value: 0.048646893

Note that $p_1 = 0.05598 > 0.05, p_2 = 0.2884 > 0.05, p_{12} = 0.0486 < 0.05$. Hence, Case (j) holds.

*Constructions for case (jj): F-test does not reject $H_0^F$, while t-tests reject $H_0^{t_1}$ or $H_0^{t_2}$ with the same significance*

Table 6A: Generated random sample for Case $(jj)$

| X1 | X2 | Y |
|---|---|---|
| 1.173699045 | 1.507797593 | 1.693611518 |
| 1.527866588 | 1.204880159 | 1.719565524 |
| -0.237756887 | 0.321525784 | 2.313343543 |
| 0.424876707 | 0.372472796 | 2.215619921 |
| 0.155008273 | -0.382097849 | 1.752313506 |
| 0.078297635 | 0.202406996 | 1.985018225 |
| -0.739378749 | -1.77490523 | 1.280511608 |
| -0.325947264 | -0.170739193 | 1.751709441 |
| 0.057639294 | 0.025498039 | 2.285726127 |
| 0.317517151 | 0.439566564 | 1.809984615 |

Table 6B: Regression output for the sample in Table 6A

| Coefficients | Estimate | Std. Error | t value | $Pr(>|t|)$ |
|---|---|---|---|---|
| (Intercept) | 1.927174047 | 0.094242869 | 20.44901714 | 1.67739E-07 |
| X1 | -0.534309107 | 0.266285323 | -2.006528569 | 0.084798865 |
| X2 | 0.478129512 | 0.20172263 | 2.37023239 | 0.049589266 |

Residual standard error: 0.271926659 on 7 degrees of freedom
Multiple R-squared: 0.445574672, Adjusted R-squared: 0.287167435560639
F-statistic: 2.812842909 on 2 and 7 DF, p-value: 0.126896814

Here $p_1 = 0.084799$, $p_2 = 0.049589$, and $p_{12} = 0.126897$. Hence $F$-test does not reject $H_0^F$, while $t$-tests reject $H_0^{t_1}$ and $H_0^{t_2}$ with significance $\alpha = 0.1$.

The next example also illustrates case $(jj)$ and, in addition, shows that $F$-test can be not inclusion monotone on the set of predictors.

Table 7A: Another generated random sample for Case ($jj$)

| X1 | X2 | Y |
|---|---|---|
| 1.568562319 | 0.927834903 | 1.612462698 |
| 1.48286001 | 1.09773946 | 2.466033052 |
| -0.573115658 | 0.981537183 | 2.518881417 |
| -0.050008016 | -0.49329821 | 1.301806858 |
| 0.165268254 | -0.397500853 | 1.436310825 |
| -0.306203404 | -0.193130393 | 2.072432714 |
| -0.399941489 | -0.096035236 | 1.573998771 |
| 0.21069356 | 0.603984432 | 1.827021014 |
| 0.431810105 | -0.383312909 | 1.933014077 |
| 0.080628207 | 0.231611299 | 2.002964703 |

Table 7B: Regression output for the sample in TAble 7A

| Coefficients | Estimate | Std. Error | t value | $Pr(>|t|)$ |
|---|---|---|---|---|
| (Intercept) | 1.8040579 | 0.1149434 | 15.69519 | 1.0317e-06 |
| X1 | -0.1814855 | 0.1720305 | -1.05496 | 0.326488 |
| X2 | 0.5168505 | 0.2013693 | 2.56668 | 0.037188 |

Residual standard error: 0.3324189 on 7 degrees of freedom
Multiple R-squared: 0.486189, Adjusted R-squared: 0.3393858
F-statistic: 3.311843 on 2 and 7 DF, p-value: 0.09723261

Note that $p_2 = 0.037188 < 0.05 < 0.097233 = p_{12}$. Hence, Case ($jj$) holds. Furthermore, eliminating predictor $X1$ yields the following SLR table:

Table 7C: Regression output for the sample in Table 7A with independent variable $X2$ only

| Coefficients | Estimate | Std. Error | t value | $Pr(>|t|)$ |
|---|---|---|---|---|
| (Intercept) | 1.7797247 | 0.1133974 | 15.69458 | 2.7117e-07 |
| X2 | 0.4157529 | 0.1783505 | 2.33110 | 0.048079 |

Residual standard error: 0.3347572 on 8 degrees of freedom
Multiple R-squared: 0.404497, Adjusted R-squared: 0.330059
F-statistic: 5.434026 on 1 and 8 DF, p-value: 0.04807907

We observe again that $p'_2 = 0.0480790 < 0.05 < 0.09723261 = p_{12}$. Thus, $F$ test states that both predictors $X1$ and $X2$ are insignificant, while $X2$ alone is significant at $\alpha = 0.05$.

## Concluding remarks

Both, ANOVA and $TK$ multiple comparisons tests with $k$ groups may result in logical contradictions when $k > 2$, even if $INAH$ assumptions hold. So, the good old approach of using pairwise comparisons instead of multiple ones is a bit slower but more reliable. Furthermore, all contradictions disappear if we replace the ANOVA and $TK$ tests by their pairwise versions, applying them for any pair of groups $i, j \in \{1, \cdots, k\}$ with $i \neq j$. Then, by Theorem 1, these two tests become equivalent.

Similar contradictions appear for the Linear Regression with the number of predictors $k > 1$ (MLR). Already for $k = 2$, with the same level of significance $\alpha$, it may happen that $t$-test rejects $H_0^{t_1} : \beta_1 = 0$, while $F$-test fails to reject the stronger null hypothesis $H_0^F : \beta_1 = \beta_2 = 0$.

In general, estimating the quality of a prediction made by ANOVA or MLR seems much more doubtful than the prediction itself.

## Acknowledgements

This paper was prepared within the framework of the HSE University Basic Research Program.